\documentclass[12pt,reqno]{amsart}
\usepackage{xy,epsfig,amsfonts,amssymb}

\catcode `\@=11

\def\numberbysection{\@addtoreset{equation}{section}
\renewcommand{\theequation}{\thesection.\arabic{equation}}}
\numberbysection
\def\subsubsection{\@startsection{subsubsection}{3}%
  \normalparindent{.5\linespacing\@plus.7\linespacing}{-.5em}%
  {\normalfont\bfseries}}
\def\cal{\mathcal}
\def\d{{\mathrm d}}
\def\F2{\mathbb F_2}

\def\E{{\mathcal E}}
\DeclareMathOperator{\Sq}{Sq}
\DeclareMathOperator{\Id}{Id}
\DeclareMathOperator{\Com}{{\mathcal C}om}
\DeclareMathOperator{\Ker}{Ker}
\DeclareMathOperator{\Ima}{Im}
\DeclareMathOperator{\Lev}{{\mathcal L}ev}
\DeclareMathOperator{\Free}{{\mathcal F}ree}
\def\LevAC{\Lev^{AC}}

\xyoption{all}
\CompileMatrices

\begin{document}
\title{Adem-Cartan operads}
\author[D. Chataur, \ M. Livernet]
{David Chataur$^*$\hfill Muriel Livernet \\
Universit\'e d'Angers \hfill Universit\'e Paris 13}

\footnotetext{$^*$ Partially supported by Marie Curie grant \#HPMF-CT-2001-01179
and european TMR ``Modern Homotopy Theory'' \#HPRN-CT-1999-00119 and
Bourse postdoctorale R\'egion pays de loire}

\date{\today}

\email{dchataur@tonton.univ-angers.fr\\
livernet@math.univ-paris13.fr}

\begin{abstract}
In this paper, we introduce  Adem-Cartan operads and prove that
the cohomology of any algebra over such an operad is an unstable
level algebra over the extended Steenrod algebra. Moreover we
prove that this cohomology is endowed with secondary cohomological
operations.\\

\noindent {\bf MSC (2000): } 55-xx; 55S10; 18D50. \\
\noindent {\bf Keywords:} cohomology operations, Adem-Cartan relations,
operads, $E_\infty$-structures, level algebras.
\end{abstract}

\maketitle

\section*{Introduction}

The Steenrod Algebra $\mathcal A_p$ is one of the main
computational tool of homotopy theory. Steenrod's operations were
first introduced by N.E. Steenrod in $1947$ \cite {St} for $p=2$
and for an odd prime in $1952$ \cite {St2}. The relations between
these cohomological operations were determined by J. Adem \cite
{Ade} and H. Cartan \cite {Ca}. Cartan's proof relies on the
computation of the singular cohomology of the Eilenberg-Mac Lane
spaces at the prime $p$. Adem's proof is based on the computation
of the homology of the symmetric group $\Sigma_{p^2}$ acting on
$p^2$ elements at the prime $p$. In the sixties, Adams introduced
secondary cohomological operations \cite {Ad}, which are an
efficient tool to deal with realisability problems of unstable
modules over $\cal A_2$. In this paper we extend these results to
a more algebraically framework, at the prime 2.

\medskip

The purpose of this paper is to describe a family of operads,
named Adem-Cartan operads, such that the cohomology of an algebra
over such an operad is an unstable level algebra over the extended
Steenrod algebra $\cal B_2$ (see corollary \ref{structureB2}). By a
{\it level algebra}, we mean a commutative algebra (not
necessarily associative) satisfying the following $4$-terms
relation:
$$(a*b)*(c*d)=(a*c)*(b*d).$$
There is an operad, denoted by $\Lev$, such that level algebras are 
algebras over the operad $\Lev$. An unstable
algebra over $\cal B_2$ is a $\cal B_2$-module satisfying the
usual conditions, namely the Cartan formula and Adem relations.
In order to describe the family of Adem-Cartan operads, the cofibrant operad
$\LevAC$ is introduced, obtained by the process of attaching cells
from the standard bar resolution of $\F2$, such
that $\LevAC\rightarrow \Lev$ is a fibration. More precisely, the
operad $\LevAC$ is generated by some 2-cells, $e_i\in
\LevAC(2)^{-i}$, corresponding to the usual $\cup_i$-products and
some 4-cells $G^m_n \in\LevAC(4)^{-n}$. An {\it Adem-Cartan operad}
is an operad $\cal O$ together with a morphism of operads
$\LevAC\rightarrow \cal O$ satisfying 
some conditions (see definition \ref{ADO}). We prove
that the cohomology of a $\cal O$-algebra is an unstable algebra
over $\cal B_2$. Namely, the cells $e_i$ are responsible for the
existence of Steenrod squares whereas the cells $G^m_n$ are
responsible for the relations between them (Adem and Cartan
relations). Moreover the cohomology of a $\LevAC$-algebra $A$ is a
$\LevAC$-algebra itself (this structure is however not natural).
Hence the 4-cells $G^m_{n+1}$ yield secondary cohomological
operations $\theta^{m,n}:H^q(A)\rightarrow H^{4q-n-1}(A)$. We
prove that given an Adem relation $\sum_i
\Sq^{m_i}\Sq^{n_i}(a)=0$, there exist maps $\psi^{m,n}$ from
$\cap_i \Sq^{n_i} \subset H^*(A)$ to $H^*(A)/\sum_i\Ima\Sq^{m_i}$
defined at the cochain level. We prove that these two maps coincide,
that is $\psi^{m,n}(a)=[\theta^{m,n}(a)]$, for
$a\in\cap_i\Sq^{n_i}$ (see theorem \ref{scogene})

\medskip

Note that we recover some classical results on topological spaces.
Since $E_\infty$-operads are Adem-Cartan operads (see
\ref{uptohomo}), any algebra over an $E_\infty$-operad 
is an unstable algebra
over the extended Steenrod algebra (see \cite{May}, \cite{KM}).
Furthermore since the cochain complex of a topological space
$C^*(X;\F2)$ is an algebra over a $E_\infty$-operad (see \cite{HS}), then it has
secondary cohomological operations. Thanks to the work of
Kristensen \cite{Kr}, we prove that there operations coincide with
Adams operations. Moreover, we can extend these operations in a
non natural way to secondary cohomogical operations on the whole cohomology
(see theorems \ref{AdamsId} and \ref{cohop2}).

\medskip

The paper is presented as follows. Section 1 contains the
background needed. In Section 2, we set the construction of the
operad $\LevAC$, define Adem-Cartan operads and prove that
$E_\infty$-operads lie in this family. Section 3 is devoted to
the fondamental theorem \ref{FOND} and its corollary \ref{structureB2}. 
Section 4
is concerned with secondary cohomological operations and section 5
is devoted to proofs of technical lemmas stated in the different
sections.

\section{Recollections}

The ground field is $\F2$. In this article, a {\it vector space}
means a differential $\mathbb Z$-graded  vector space over $\F2$,
where the differential is of degree 1.

The symbol $\Sigma_n$ denotes the symmetric group acting on $n$
elements. Any $\sigma\in \Sigma_n$ is written
$(\sigma(1)\ldots\sigma(n))$.

\subsection{Operads.}(\cite {GJ}, \cite {GK}, \cite {KM}, \cite {Lo})
A (right) {\it $\Sigma_n$-module} is a (right)
$\mathbb F_2[\Sigma_n]$-differential
graded module. A {\it $\Sigma$-module} $\mathcal M=\{\mathcal M(n)\}_{n>0}$
is a collection of $\Sigma_n$-modules. Any $\Sigma_n$-module $M$
gives rise to a $\Sigma$-module $\mathcal M$ by setting $\mathcal
M(q)=0$ if $q\not =n$ and $\mathcal M(n)=M$.

\bigskip

An {\it operad} is a right $\Sigma$-module $\{\mathcal
O(n)\}_{n>0}$ such that $\mathcal O(1)=\mathbb F_2$, together with
composition products:

\begin{eqnarray*}
\mathcal O(n)\otimes\mathcal O(i_1)\otimes\ldots\otimes\mathcal
O(i_n)&\longrightarrow &\mathcal O(i_1+\ldots +i_n) \\
o\otimes o_1\otimes\ldots\otimes o_n &\mapsto&
o(o_1,\ldots,o_n).
\end{eqnarray*}

These compositions are subject to associativity conditions,
unitary conditions and equivariance conditions with respect to the
action of the symmetric group. The equivariance conditions write:

\begin{align*}
o(a_1\cdot\tau_1,\ldots,a_n\cdot\tau_n)&=&
o(a_1,\ldots,a_n)\cdot(\tau_1\times\ldots\times\tau_n) \\
(o\cdot\sigma)(a_1,\ldots,a_n)&=&o(a_{\sigma^{-1}(1)},\ldots,a_{\sigma^{-1}(n)})\cdot\sigma(i_1,\ldots,i_n),
\end{align*}
where $\tau_1\times\ldots\times\tau_n$ is the permutation of
$\Sigma_{i_i+\ldots+i_n}$ such that $\tau_1$ operates on the first
$i_1$ terms, $\tau_2$ on the next $i_2$ terms and so on; the permutation
$\sigma(i_1,\ldots,i_n)$ in
$\Sigma_{i_i+\ldots+i_n}$ operates as $\sigma$ on the $n$-blocks
of size $i_k$.

For instance, for any $\sigma,\mu,\nu\in \cal O(2)$, one has
\begin{equation}\label{flip1}
(\sigma\cdot (21))(\mu,\nu)=\sigma(\nu,\mu)\cdot (3412)
\end{equation}

There is another definition of
operads via $\circ_i$ operations
$$\circ_i : \mathcal O(n)\otimes \mathcal O(m)\longrightarrow \mathcal
O(n+m-1),$$ where $p\circ_i q$ is $p$ composed with $n-1$ copies of the
unit $1\in\mathcal O(1)$ and with $q$ at the $i$-th position.

\medskip

The forgetful functor from the category of operads to the category of
$\Sigma$-modules has a left adjoint : the {\it free operad functor}, denoted
by
$\Free$.

\bigskip

An {\it algebra over an operad $\mathcal O$} or a {\it $\mathcal O$-algebra} $A$
is a vector space together with evaluation maps
\begin{eqnarray*}
\mathcal O(n)\otimes A^{\otimes n}&\longrightarrow &A \\
o\otimes a_1\otimes\ldots\otimes a_n &\mapsto& o(a_1,\ldots,a_n)
\end{eqnarray*}
These evaluation maps are subject to associativity conditions and
\hfill\break equivariance conditions. These equivariance
conditions write:

\begin{align*}
(o\cdot\sigma)(a_1,\ldots,a_n)=o(a_{\sigma^{-1}(1)},\ldots,a_{\sigma^{-1}(n)}).
\end{align*}

\smallskip

A {\it  graded algebra over an operad $\mathcal O$} is a $\mathcal
O$-algebra whose differential is null.

\subsection{Operadic cellular attachment.}(\cite {BM}, \cite {Hi})
\label{cellattach} The category of operads is a closed model
category. Weak equivalences are qua\-si-iso\-mor\-phisms (i.e.
isomorphisms in cohomology of the underlying vector spaces) and
fibrations are epimorphisms. Cofibrations can be defined by the
left lifting property with respect to the acyclic fibrations. For
background material on closed model categories we refer to \cite
{DS}, \cite {Ho} and \cite {Qu}.
\bigskip
Any morphism of operads
$$\mathcal P\longrightarrow \mathcal Q$$
can be factorized by a
cofibration ($\xymatrix{ \ar@{>->}[r] & }$) followed by an acyclic fibration
($\xymatrix{\ar@{>>}[r]^{\sim}& }$). This factorization can
be realized using the inductive process of attaching cells (the
category of operads is cofibrantly generated \cite{Hi}).
\smallskip
An operad is {\it cofibrant} if the morphism from the initial object
$\Free(0)$ to
the operad is a cofibration. In order to produce a {\it cofibrant replacement}
to an operad $\mathcal O$, one applies the inductive process
of attaching cells to the canonical morphism $\Free(0)\longrightarrow \mathcal O
$.
\medskip

Let $S_p^n$ be the free $\Sigma_p$-module generated by $\delta t$
in degree $n$ considered as a $\Sigma$-module. Let $D_p^{n-1}$ be
the free $\Sigma_p$-module generated by $t$ in degree $n-1$ and
$dt$ in degree $n$, the differential sending $t$ to $dt$. We have
a canonical inclusion $i_{n}:S_p^{n}\longrightarrow D_p^{n-1}$ of
$\Sigma$-modules (sending $\delta t$ to $dt$). Let
$f:S^n_p\longrightarrow \mathcal O$ be a morphism of
$\Sigma$-modules. The cell $D_p^{n-1}$ is attached to $\mathcal O$
along the morphism $f$ via the following push-out:

$$\xymatrix{
{\Free}(S_p^{n}) \ar[d]_{\Free(f)} \ar[r]^{i_n} & \Free(D_p^{n-1}) \ar[d]  \\
{\mathcal O}\ \  \ar[r]^>>>>>>{i}& {\mathcal O}\coprod_{\tau}\Free(S_p^{n-1}).}
$$

The main point of this process is that $f(\delta t)$, which was a
cycle in $\cal O$,
becomes a boundary in $\cal O\coprod_{\tau}\Free(S_p^{n-1})$.

\smallskip

By iterating this process of cellular attachment, one gets a {\it
quasi-free extension}: $\xymatrix{ {\mathcal O}\ \ \ar[r]^>>>>{i}&
{\mathcal O}\coprod_{\tau}\Free(V),} $ that is if we forget the
differential on $V$ then $\mathcal O\coprod_{\tau}\Free(V)$ is the
coproduct of $\mathcal O$ by a free operad over a free graded
$\Sigma$-module $V$. Note that any cofibration is a retract of a quasi-free 
extension. A {\it quasi-free operad} is an operad which is free over a free 
$\Sigma$-module if we forget the differential.

\medskip

\noindent The following proposition will be fundamental for our applications.

\subsubsection{Proposition} \label{propcellattach}\it Let $V$ be a free graded
$\Sigma_p$-module
together with
$$d_V:V\longrightarrow \mathcal O(p)\oplus V$$
such that $d_V+d_{\mathcal O}$ is of square zero. Then
if $V$ is bounded above, the morphism
$\mathcal O\rightarrow \mathcal O\coprod_{\tau} \Free(V)$
is a cofibration. \rm

\medskip

\noindent{\bf Proof.} The boundness assumption is needed in order to build a map

$$\xymatrix{
{\mathcal O}\ \ar@{>->}[r] &{\mathcal O}\coprod_{\tau} \Free(V)}$$
by induction on the degree of $V$, using the cellular attachment process.
In order to add cells the assumption  $(d_V+d_{\mathcal O})^2=0$ is needed.~~~~\hfill{$\square$}

\subsection{Homotopy invariance principle.}
\label{invariance} Let $\mathcal O$ be a cofibrant operad. The
category of $\mathcal O$-algebras is also a closed model category
where weak equivalences are quasi-isomorphisms and fibrations are
epimorphisms \cite{BM} and \cite{Hi}.

Recall that the category of vector spaces is a closed model category, 
where weak  equivalences are
quasi-isomorphisms and fibrations are epimorphisms. In this
category all objects are fibrant and cofibrant.

The following theorem is stated in  \cite{Ch}; its proof relies
on a general result of Berger and Moerdijk about transfer of
algebraic structure  in closed model categories \cite{BM}. 
Such a result was proved in caracteristic zero by Markl \cite{Mar}
and for topological spaces by Boardman and Vogt \cite{BV}.

\subsubsection{Theorem}({\bf Homotopy invariance principle.}) \cite{Ch}

{\it Let $\mathcal  O$ be a cofibrant operad and assume that the
morphism of vector spaces
$$f:X\longrightarrow Y$$
is a weak equivalence between vector spaces. Assume that $X$ is a
$\mathcal O$-algebra. Then $Y$ is also provided with a $\mathcal
O$-algebra structure. For any cofibrant replacement $\tilde X$ of
$X$, there exists a sequence of quasi-isomorphisms of $\mathcal
O$-algebras:
$$X\longleftarrow \tilde X\longrightarrow Y.$$}
such that the diagram commute in cohomology:
$$\xymatrix{
& H^*(\tilde X) \ar[dl]\ar[dr]& \\
H^*(X)  \ar[rr]^{H^*(f)} & & H^*(Y)}$$

\subsubsection{Corollary.} \label{corinvariance}\it Let $\cal O$ be a cofibrant
operad.  Let $C$ be a $\mathcal O$-algebra and let $H$ be its
cohomology. Then $H$ is a $\mathcal O$-algebra and  there is a
sequence of quasi-iso\-mor\-phisms of $\mathcal O$-algebras
$$\xymatrix{
& \tilde C \ar[dl]_{\psi}\ar[dr]^{\Phi}& \\
C   & & H^*(C)}$$ where $\tilde C$ is a cofibrant replacement of
$C$ and such that $H^*(\Phi)=H^*(\psi)$. \rm

\medskip

\noindent{\bf Proof.} Set $X=C, Y=H^*(C)$. Choose a decomposition
of vector spaces of $X=Z\oplus Q$, where $Z$ is the kernel of the
differential and let $f:X\rightarrow Y$ be the composite of the
projections onto $Z$ and onto $H$. Then $f$ is a weak equivalence
of vector spaces, since $H^*(f)=\Id$, and we can apply the
previous theorem. ~~~\hfill{$\square$}

\section{Adem-Cartan operads}

In this section, the operad $\LevAC$ is built, 
as a first step towards a resolution of the operad $\Lev$. 
The latter operad governs level algebras which are commutative algebras satisfying
the 4-terms relation (\ref{rlev}). The idea of introducing level algebras
instead of commutative algebras in order to deal with Cartan and Adem 
relations is inspired by the fact that these relations are not conditioned
by the associativity of the product. In this section we define also
the family of Adem-Cartan operads and prove that
$E_\infty$-operads belong to this family.

\subsection{The operad $\Lev$}\label{Lev}  A {\it level algebra} A is a vector
space together with a commutative product $*$ (non necessarily associative) satisfying
the relation
\begin{equation} \label{rlev}
(a*b)*(c*d)=(a*c)*(b*d), \ \forall a,b,c,d\in A.
\end{equation}

\subsubsection{Definition} Let $\F2$ be the trivial representation of $\Sigma_2$
(generated by the operation $\mu$) and $R_{\Lev}$ be the
sub-$\Sigma_4$-module of
$\Free(\F2)(4)$ generated by the elements $\mu(\mu,\mu)\cdot (\Id+\sigma)$
for all $\sigma\in\Sigma_4$. Then the operad $\Lev$ is the operad
$$\Lev=\Free(\F2)/<R_{\Lev}>.$$
Algebras over this operad are level algebras.

\subsubsection{Remark} \label{comlev} Since a commutative and associative
algebra is trivially a level algebra,
there is a morphism of operads
$$\Lev\longrightarrow \Com,$$
where $\Com$ denotes the operad defining commutative and associative algebras.

\subsubsection{Definition} \label{shape} For any $\Sigma_2$-module M, the
vector space \hfill\break
$\Free(M)(4)$ is a direct sum of two $\Sigma_4$-modules: the one
indexed by trees of {\it shape 1} denoted by $\mathbb T_1$, that is the
$\Sigma_4$-module generated by all the compositions
$\mu\circ_1\gamma$ for $\mu\in M$ and
$\gamma\in\Free(M)(3)$ and the one indexed
by trees of {\it shape 2} denoted by $\mathbb T_2$, that is the
$\Sigma_4$-module generated by all the compositions
$\mu(\nu,\eta)$ for $\mu,\nu,\eta \in M$.

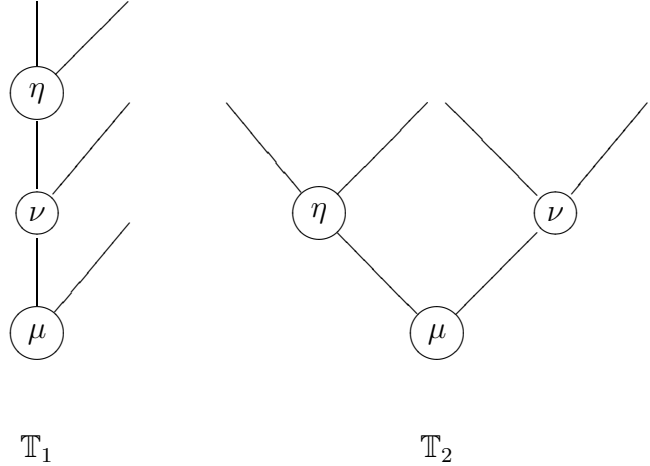
$$\vcenter{\xymatrix{
&&&&& \\
*++[o][F-]{\eta} \ar@{-}[ur]
 \ar@{-}[u] \ar@{-}[d] &&&&&&    \\
*++[o][F-]{\nu} \ar@{-}[ur]\ar@{-}[d]& &  & *++[o][F-]{\eta}
\ar@{-}[ul]
 \ar@{-}[ur] \ar@{-}[dr] &&*++[o][F-]{\nu} \ar@{-}[ur]
 \ar@{-}[ul] \ar@{-}[dl] &     \\
*++[o][F-]{\mu}\ar@{-}[ur]& & & & *++[o][F-]{\mu}& & \\
{\mathbb T_1} & & &&   {\mathbb T_2}  }}
$$

As an example, since there is only one ge\-ne\-ra\-tor $\mu \in
\Free(\F2)(2)$, the dimension of $\mathbb T_1(\Free(\F2))$ is 12
and the dimension of $\mathbb T_2(\Free(\F2))$ is 3,
whereas the dimension of $\mathbb T_1(\Lev)$ is 12
and the dimension of $\mathbb T_2(\Lev)$ is 1.

\subsection{Construction of $\Lev^{AC}$.}\label{constructionlevAC}
The construction of this operad is done by the process of 
attaching cells (see \ref{cellattach}) to the
standard bar resolution $\E$ of ${\mathbb F}_2[\Sigma_2]$.

\subsubsection{The standard bar resolution} \label{barres} Let $(21)$
be the non
trivial permutation of $\Sigma_2$.
The standard bar re\-so\-lu\-tion of $\Sigma_2$ over $\mathbb F_2$
is given by
\begin{align*}
\E^{-i}&=\begin{cases}
<e_i,e_i\cdot (21)>&\mathrm{if \ } i\geq 0, \\
0, & \mathrm{if\ } i<0,
\end{cases} \\
\mbox{and}\  d(e_i)&=e_{i-1}+ e_{i-1}\cdot (21), \mathrm{with\ } e_{-1}=0.
\end{align*}

Denote again by $\E$ the free operad generated by
the $\Sigma$-module $\E$. 

\medskip

Since $\Lev=\Free(\F2)/<R_{\Lev>}$, there is
a fibration of operads:
$$p:\xymatrix{\E \ar@{>>}[r]^>>>>{\sim}&\Free (\F2) \ar@{>>}[r] &\Lev.}$$
By definition of $\Lev$, $p(n)$ is clearly a quasi-isomorphism for
$n<4$. Furthermore $p(4)$ is a quasi-isomorphism on pieces of
shape 1 (see \ref{shape}). 

\medskip

Note that for any operad $\cal O$ such that $\cal O(2)$ is 
a $\Sigma_2$-projective
re\-so\-lu\-tion of $\F2$,  there exists a
morphism of operads $m: \E \rightarrow \cal O$ such that $m(2)$ is a homotopy equivalence.
Then the image $m(e_i)\in \cal O(2)$ is non zero. In the sequel, this image will be denoted
also by $e_i$.

\subsubsection{Notation}--We use May's convention: for any integers
$i$ and $j$ the symbol
$(i,j)$ denotes $\frac{(i+j)!}{i!j!} \in \F2$, if $i\geq 0$ and
$j\geq 0$ and $(i,j)=0$ otherwise. If the 2-adic expansion of $i$
and $j$ is $i=\sum a_k 2^k$ and $j=\sum b_k 2^k$ then $(i,j)=0 \in
\F2$ if and only if there exists $k$ such that $a_k=b_k=1$.

\medskip

For any $\sigma, \tau, \nu \in \E(2)$ we denote $\sigma\cdot
(21)(\tau,\nu)$ by $\sigma(\tau,\nu)\cdot F$ which is equal the to
$\sigma(\nu,\tau)\cdot (3412)$ (see (\ref{flip1})).

\subsubsection{Definition} \label{defsuite}
\begin{itemize}
\item[a)] Let define some elements
$[u_n^m]_x\in\E(4)^{-n-x}$, for $m>0$ which are a combination of elements of shape 2. For  any $m$ such that $2^k\leq m\leq 2^{k+1}-1$,
\begin{equation}\label{defun}
\begin{split}
[u^m_0]_x&= e_x\cdot (21)^{m-1}(e_0,e_0),\ {\rm and\ for\ } n>0,\\
[u^m_n]_x&=\sum_{i=0}^{2^{k+1}-1} \sum_{ 0\leq 2^{k+1}\delta-i\leq n}(n-m+i,m-1)[ \\
&(i,m) e_x\cdot (21)^{m-1} (e_{2^{k+1}\delta-i},e_{n+i-2^{k+1}\delta}) + \\
&(i-1,m) e_x\cdot (21)^{m-1}(e_{2^{k+1}\delta-i},e_{n+i-2^{k+1}\delta}\cdot(21))]
\end{split}
\end{equation}

\item[b)]  Let define some elements $\alpha_{n,p}\in \E(4)^{-n}$, for
$p\in{\mathbb Z}$ 
\begin{equation}\label{defalpha}
\begin{split}
\alpha_{n,p}&=\sum_{s=0}^{n-p} [u_{s+p}^{s+1}]_{n-s-p}+
\sum_{s=0}^{n-p-1} [u_{s+p+1}^{s+1}]_{n-s-p-1}\cdot (3412) \\
\end{split}
\end{equation}
\end{itemize}

\subsubsection{Proposition}\label{suite}\it
The $[u_n^m]_x$'s satisfy the following properties:

\begin{align}
[u_n^m]_x=&0, \ ~{\rm for}~ 0<n<m  \label{unmzero}\\
[u_m^m]_x=&e_x\cdot (21)^{m-1}(e_0,e_m)+e_x\cdot (21)^{m-1}(e_m,e_0\cdot(21))
\label{umm}\\
\begin{split}
\d [u^{m+1}_{n+1}]_x=&[u^{m+1}_{n}]_x(\Id+(2143))+[u^{m}_{n}]_x(F+(4321))\\
&+[u^{m+1}_{n+1}]_{x-1}(\Id+F)
\end{split} \label{du}
\end{align}

The $\alpha_{n,p}$'s satisfy the following property:
$$\alpha_{n,p}=0,\ \mbox{if}\ p<0,$$
\begin{align}
\d \alpha_{n,p}=&\alpha_{n-1,p-1}\cdot (\Id+(2143))+\alpha_{n-1,p}\cdot (\Id+(4321))
\label{dalpha}
\end{align}

\rm

The proof of this proposition is postponed to the last section.

\subsubsection{Examples} \label{u1} By definition we have for all
$n\geq 0$
\begin{align*}
[u_n^1]_x=&\sum_k e_x(e_{2k},e_{n-2k})+
\sum_l e_x(e_{2l+1},e_{n-2l-1}\cdot (21))\\
=&e_x(\psi(e_n))\\
\end{align*}
where $\psi$ is the coproduct in the standard bar resolution $\E$ of $\F2$
(see \cite {May}, \cite {BF}). And for $n\geq 0$, one has
\begin{align*}
[u_n^2]_x=\begin{cases}
\sum_{ 0\leq 2\delta\leq n} e_x\cdot (21) 
(e_{2\delta},e_{n-2\delta}\cdot (21)^{\delta}) &
{\rm if\  n\  even}\ \\
\sum_{ 0\leq 4\delta-1\leq n}
e_x\cdot (21) (e_{4\delta-1},e_{n+1-4\delta}((12)+(21))\ & 
{\rm if\ n\  odd}
\end{cases}
\end{align*}

\subsubsection{Theorem}\label{cell} \it There exists a cofibrant
operad
$\LevAC$ satisfying the following properties:
\begin{itemize}
\item[a)] $\LevAC(2)=\E$;
\item[b)] there is a fibration $f:\LevAC\rightarrow \Lev$ such that
$f(n)$ is a quasi-isomorphism for $n<4$ and induces an isomorphism \\
$H^0(\LevAC(n))\simeq \Lev(n)$;
\item[c)] there are elements $G^m_n\in\LevAC(4)$ of degree $-n$,  satisfying
\begin{align}
\begin{split}
G^m_n&=0 \ {\it for} \ m\leq 0 \ {\it or}\ n<m \\
\d G^{m}_{m}&=G^{m-1}_{m-1}\cdot (3214)(\Id+(2143))+
G^{m-2}_{m-1}\cdot (\Id+(4321))  \\
&+ \alpha_{m-1,m-1-p}+\alpha_{m-1,p}\cdot (3214) \\
&\qquad {\it and\ for\ } n>m \\
\d G^{m}_{n}&=G^{m}_{n-1}\cdot (\Id+(2143))+
G^{m-2}_{n-1}\cdot (\Id+(4321))  \\
&+\alpha_{n-1,n-1-p}+\alpha_{n-1,p}\cdot (3214)
\end{split}\label{cellmn}
\end{align}
where $p$ is the integer part of $\frac{m-1}{2}$.

\end{itemize}\rm

\medskip

\noindent {\bf Proof.} The proof consists of building a sequence of cofibrant
operads $\LevAC_m$, $m\geq 0$, satisfying a) and b) where
$\LevAC_m\rightarrow \LevAC_{m+1}$ is a cofibration obtained by operadic
cellular attachment.
More precisely $\LevAC_0=\E$; assume $\LevAC_m$ is built and let $V^{m+1}$ be the
free graded $\Sigma_4$-module generated by elements $G_n^{m+1}$ of degree $-n$ for $n\geq m+1$
with $d:V^{m+1}\rightarrow \LevAC_m(4)\oplus V^{m+1}$ defined by the relation (\ref{cellmn}).
If $d^2=0$, according to proposition \ref{propcellattach}, then
$$\LevAC_m\rightarrow \LevAC_m
\coprod_{\tau}\Free(V^{m+1})=: \LevAC_{m+1}$$
is a cofibration. The main ingredient to prove that $d^2=0$ is the
relation (\ref{dalpha}).

\smallskip

\noindent{\it First case: $n>m+1$}; then $n-1>m$ and $n-1>m-2$
\begin{align*}
(\d)^2(G^m_{n})&=\d (G^{m}_{n-1}\cdot (\Id+(2143))+
G^{m-2}_{n-1}\cdot (\Id+(4321)))\\
&+\d(\alpha_{n-1,n-1-p}+\alpha_{n-1,p}\cdot (3214)) \\
&=G^m_{n-2}\cdot (\Id+(2143))^2+ G^{m-2}_{n-2}\cdot(\Id+(4321))(\Id+(2143))\\
&+G^{m-2}_{n-2}\cdot(\Id+(2143))(\Id+(4321))+G^{m-4}_{n-2}\cdot (\Id+(4321))^2 \\
&+\alpha_{n-2,n-2-p}\cdot(\Id+(2143))+\alpha_{n-2,p}\cdot(3214)(\Id+(2143))\\
&+\alpha_{n-2,n-2-(p-1)}\cdot(\Id+(4321))\\
&+\alpha_{n-2,p-1}\cdot(3214)(\Id+(4321))\\
&+\d(\alpha_{n-1,n-1-p}+\alpha_{n-1,p}\cdot (3214)) \\
\end{align*}
Since $(2143)$ and $(4321)$ commute and are permutations of order 2
the 2 first lines vanish. Furthermore
$$\d \alpha_{n-1,n-1-p}=\alpha_{n-2,n-2-p}\cdot (\Id+(2143))+
\alpha_{n-2,n-1-p}\cdot(\Id+(4321)),$$
and $(3214)(\Id+(4321))=(\Id+(2143))(3214)$ implies
\begin{multline*}
\d \alpha_{n-1,p}\cdot (3214)=\alpha_{n-2,p-1}\cdot(3214)(\Id+(4321))\\
+\alpha_{n-2,p}\cdot(3214)(\Id+(2143))
\end{multline*}
which yields  $(\d)^2(G^m_n)=0$ for $n>m+1$.

\smallskip

\noindent{\it Second case: $n=m+1$}, then $n-1=m$ and $n-1>m-2$; thus terms in $\alpha_{u,v}$
will vanish as before and
\begin{align*}
(\d)^2(G^m_{m+1})=&G^{m-1}_{m-1}\cdot(3214) (\Id+(2143))^2+ G^{m-2}_{m-1}\cdot(\Id+(4321))(\Id+(2143))\\
&+G^{m-2}_{m-1}\cdot(\Id+(2143))(\Id+(4321))+G^{m-4}_{m-1}\cdot (\Id+(4321))^2 \\
=&0 \\
\end{align*}

\noindent{\it Last case: $n=m$}, then
\begin{align*}
(\d)^2(G^m_{m})=&\d (G^{m-1}_{m-1}\cdot(3214) (\Id+(2143))+
G^{m-2}_{m-1}\cdot (\Id+(4321)))\\
&+\d(\alpha_{m-1,m-1-p}+\alpha_{m-1,p}\cdot (3214)) \\
=&G^{m-2}_{m-2}\cdot [(3214)(\Id+(2143))]^2+ G^{m-3}_{m-2}\cdot(\Id+(4321))(3214)(\Id+(2143))\\
&+G^{m-2}_{m-2}\cdot(\Id+(2143))(\Id+(4321))+G^{m-4}_{m-2}\cdot (\Id+(4321))^2 \\
&+\alpha_{m-2,m-2-p'}\cdot(3214)(\Id+(2143))+\alpha_{m-2,p'}\cdot(3214)^2(\Id+(2143))\\
&+\alpha_{m-2,m-2-(p-1)}\cdot(\Id+(4321))+\alpha_{m-2,p-1}\cdot(3214)(\Id+(4321))\\
&+\d(\alpha_{m-1,m-1-p}+\alpha_{m-1,p}\cdot (3214)) \\
\end{align*}
where $p'$ is the integer part of $\frac{m-2}{2}$ and $p$ the integer
part of $\frac{m-1}{2}$.

It's easy to check that the 2 first lines vanish. Furthermore, if $m=2k+1$ is odd
then $p=k$, $p'=k-1$ and the last lines write
\begin{align*}
&+\alpha_{2k-1,k}(\Id+(4321))\cdot(3214)+\alpha_{2k-1,k-1}(\Id+(2143))\\
&+\alpha_{2k-1,k}\cdot(\Id+(4321))+\alpha_{2k-1,k-1}\cdot(\Id+(2143))(3214)\\
&+\d(\alpha_{2k,k}+\alpha_{2k,k}\cdot (3214))=0, \\
\end{align*}
and if $m=2k+2$ is even, then $p=k$, $p'=k$ and the last lines write
\begin{align*}
&+\alpha_{2k,k}(\Id+(4321))\cdot(3214)+\alpha_{2k,k}(\Id+(2143))\\
&+\alpha_{2k,k+1}\cdot(\Id+(4321))+\alpha_{2k,k-1}\cdot(\Id+(2143))(3214)\\
&+\d(\alpha_{2k+1,k+1}+\alpha_{2k+1,k}\cdot (3214))=0 \\
\end{align*}

\smallskip

Let us prove that the operads $\LevAC_m$ satisfy b). According to
\ref{barres} $\E(n)\rightarrow \Lev(n)$ is a quasi-isomorphism for
$n<4$, so is   $\LevAC_m(n)\rightarrow \Lev(n)$ for all $m$.

\smallskip

To prove that $f_m:\LevAC_m\rightarrow \Lev$ induces an
isomorphism $H^0(f_m)$ for all $m>0$ it suffices to prove that
$H^0(f_1(4))$ is an isomorphism; if so, since the relations defining
$\Lev$ are generated by a $\Sigma_4$-module, then
$H^0(\LevAC_1(n)) \rightarrow \Lev(n)$  is an isomorphism, for all
$n$; if $m>1$, we do not introduce cells in degree $-1$, hence
$H^0(f_m(n))$ is an isomorphism for all $n$.

\begin{multline*}
H^0(\LevAC_1(4))={\mathbb T}_1(\E^0(4))/\d({\mathbb T}_1(\E^1(4)))
\bigoplus \\
 {\mathbb T}_2(\E^0(4))/\d[{\mathbb T}_2(\E^1(4))\oplus
G_1^1\cdot\F2(\Sigma_4)]
\end{multline*}
the first summand being isomorphic to ${\mathbb T_1}(\Lev(4))$ (see \ref{barres}). To prove that the second summand is isomorphic to ${\mathbb T_2}(\Lev(4))$
it suffices to prove that it is 1-dimensional (see \ref{shape}). Let
$X=e_0(e_0,e_0)\in \E^0(4)$ and $\overline X$ its class in the second summand.
\begin{itemize}
\item[i)] For all $\sigma$ in the dihedral group $D_1$ ($\sigma$
satisfies $\{\sigma(1),\sigma(2)\} \subset \{1,2\}\cup \{3,4\}$),
there exists $u=e_i(e_j,e_k), i+j+k=1,$ such that $\d u=X+X\cdot
\sigma$. Hence for all $\sigma\in D_1$ we have $\overline{X\cdot
\sigma}=\overline X$. \item[ii)] Since $\d G^1_1=X+X\cdot (3214)$,
for all $\sigma\in D_1$ we have $\overline{X\cdot
(3214)\sigma}=\overline X$, i.e. for all $\tau\in D_2=\{\tau
|\{\tau(1),\tau(2)\} \subset \{1,4\}\cup \{2,3\}\}$ we have
$\overline{X\cdot\tau}=\overline X$. \item[iii)]
 Finally equality $(3124)=(2134)(3214)$ implies $\overline{X\cdot (3124)\sigma}=\overline X$
for all $\sigma\in D_1$, or  for all $\tau\in D_3=\{\tau |\{\tau(1),\tau(2)\}
\subset \{1,3\}\cup \{2,4\}\}$, $\overline{X\cdot\tau}=\overline X$.
\end{itemize}
Since
$\Sigma_4=D_1\cup D_2\cup D_3$ we get the result.

Finally we define $\LevAC$ to be the limit over $m$ of $\LevAC_m$.
Then $\LevAC$ satisfies a), b) and c).~~~\hfill{$\square$}

\subsubsection{Definition} An {\it Adem-Cartan algebra} is an algebra over $\LevAC$.

\subsection{Adem-Cartan operads}

\subsubsection{Definition} \label{ADO} An {\it Adem-Cartan} operad
$\cal O$ is an operad such that  $\cal O(2)$ is a $\Sigma_2$-projective resolution of
$\F2$ together  with a
morphism of operads $m:\LevAC\rightarrow \cal O$ 
such that $m(2)$ is a quasi-isomorphism.
More precisely, there exists non trivial $e_i\in\cal O(2)^{-i}$  and 
$G_n^m \in \cal O(4)^{-n}$ satisfying the
relations (\ref{cellmn}).

\subsubsection{$E_{\infty}$-operads}\label{uptohomo}
According to remark \ref{comlev}, there is a morphism
$\Lev\longrightarrow \Com$. An $E_{\infty}$-operad  $E$
is a $\Sigma_*$-projective resolution of the operad $\Com$: for every r, 
$E(r)$
is $\Sigma_r$-projective and there exists an acyclic fibration
 $E \rightarrow \Com$. 
The structure of closed model category on operads implies that 
for any cofibrant operad
$\cal O\rightarrow \Com$ there exists a morphism $\cal O\rightarrow E$
such
that the following diagram commutes:
$$\xymatrix{
\Free(0) \ar@{>->}[d] \ar[r] &   E \ar@{>>}[d]^{\sim} \\
{\mathcal O} \ar[r] \ar[ur]& \Com. }$$ For instance there exists a
morphism $\LevAC\rightarrow E$ such that the previous
diagram commutes. Hence any $E_\infty$-operad is an Adem-Cartan
operad. In particular the algebraic
Barratt-Eccles operad $\mathcal{BE}$ (studied in \cite {BF}) is an
Adem-Cartan operad.

\section{Adem-Cartan algebras and the extended Steenrod algebra}

The aim of this section is to prove that the cohomology of an
algebra over an Adem-Cartan operad carries an action of the
extended Steenrod algebra.

\subsection{The extended Steenrod algebra $\mathcal B_2$.}

\subsubsection{Generalized Steenrod powers} (\cite {May}, \cite {Ma}) 
The {\it extended Steenrod algebra},
denoted by $\mathcal B_2$, is a graded associative algebra over
$\mathbb F_2$ generated by the generalized Steenrod squares
$\Sq^i$ of degree $i\in \mathbb Z$. These generators satisfy the Adem relation,
if $t<2s$:
$$\Sq^t \Sq^s=\sum_{i}
\binom{s-i-1}{t-2i} \Sq^{s+t-i} \Sq^i.$$
Note that negative Steenrod squares are allowed
and that $\Sq^0=\Id$ is not assumed. If $\mathcal A_2$ denotes the
classical Steenrod algebra, then
$$\mathcal A_2\cong\frac{\mathcal B_2}{<\Sq^0+\Id>}.$$

\subsubsection{Definition} As in the classical case,
an {\it unstable module over $\mathcal B_2$}, is a graded $\mathcal B_2$-module
together
with the instability condition.
\begin{align*}
\Sq^n(x)&=0 \ {\mathrm{if}}\   |x|<n
\end{align*}
{\it An unstable level algebra over $\mathcal B_2$} is a graded
level algebra $(A,*)$ which is an unstable module over $\mathcal B_2$, such that

$$\begin{cases}
\Sq^{|x|}(x)=x*x, &  \\
 \Sq^s(x*y)=\sum  \Sq^t(x)* \Sq^{s-t}(y) & ~{\rm (Cartan\  relation)}
\end{cases}$$

\subsubsection{Remark.} The category of unstable algebras over $\mathcal A_2$
is a full subcategory of the category of unstable level
algebras over $\mathcal B_2$.

\subsection{Cup-i products.} Let $A$ be an algebra over an operad
$\cal O$ such that $\cal O(2)$ is a $\Sigma_2$-projective resolution of $\F2$. 
The evaluation map
$$\cal O(2)\otimes A^{\otimes 2}\longrightarrow A,$$
defines cup-i products $a\cup_i b=e_i(a,b)$ ($e_i$'s were defined in \ref{barres}). 
Steenrod squares are defined by $\Sq^r(a)=a\cup_{|a|-r}a=e_{|a|-r}(a,a)$. Following P.
May \cite {May}, define $D_n(a)=e_n(a,a)=\Sq^{|a|-n}(a)$. In that
terminology, Adem relations read
\begin{align} \label{Ademrel2}
\sum_k (k,v-2k) D_{w-v+2k} D_{v-k}(a)&=\sum_l (l,w-2l)
D_{v-w+2l}D_{w-l}(a),
\end{align}
and Cartan relation read
\begin{align} \label{Cartanrel2}
D_n(x*y)&=\sum_{k=0}^n D_k(x)*D_{n-k}(y).
\end{align}
Note that at this stage, we do not know that an $\cal O$-algebra
satisfies these two relations.

\subsection{Main theorem}

\subsubsection{Lemma} \label{boundary}{\it Let $A$ be a graded 
algebra ($d_A=0$)
over an operad $\mathcal O$. Assume that $o\in \mathcal O(n)$ is a
boundary. Then $o(a_1,\ldots,a_n)=0$, for all $a_1,\ldots, a_n$ in
$A$.}

\medskip

\noindent {\bf Proof.} There exists $\omega\in \mathcal O(n)$ such that $d\omega=o$.
The Leibniz rule implies that
$$d_A(\omega(a_1,\ldots,a_n))=o(a_1,\ldots,a_n)+\sum_{i}
\omega(a_1,\ldots,d_Aa_i,\ldots,a_n)$$
and the result follows because $d_A=0$.~~~~\hfill{$\square$}

\subsubsection{Lemma} \label{lemmeuD} \it
Let $\cal O$ be an operad such that $\cal O(2)$ is a $\Sigma_2$-projetive
resolution of $\F2$.
For every graded  $\cal O$-algebra $A$ and for every $a\in A$,
$$[u_n^m]_x(a,a,a,a)=(n-2m,2m-1)D_xD_{\frac{n}{2}}(a),$$
where $[u_n^m]_x$ denotes the image of $[u_n^m]_x \in \E$ by
$\E\rightarrow\cal O$. \rm

\medskip

The proof is postponed to the last section.

\subsubsection{Proposition}\label{CAisLev} {\it Let $A$ be a graded Adem-Cartan
algebra, then it is a level algebra.}

\medskip

\noindent{\bf Proof.} For any operad $\cal O$ and any $\cal
O$-algebra A, $H^*(A)$ is a $H^*(\cal O)$-algebra hence a
$H^0(\cal O)$-algebra. Since $H^0(\LevAC)=\Lev$ and $H^*(A)=A$,
one gets the result.~~~~\hfill{$\square$}

\subsubsection{Theorem} \label{FOND}\it Let $A$ be a graded
Adem-Cartan algebra then $A$ is an unstable level algebra over the
extented Steenrod algebra.\rm

\medskip

\noindent {\bf Proof}. We have already proved in proposition
\ref{CAisLev} that a graded Adem-Cartan algebra is a level
algebra. Assume first that it is a module over $\cal B_2$ (Adem
relation (\ref{Ademrel2})). Hence the instability condition reads
$$\Sq^n(x)=x\cup_{|x|-n} x=0, ~{\rm if}~ |x|-n < 0.$$
The next equality is also immediate
$$\Sq^{|x|}(x)=x\cup_0 x=e_0(x,x)=x*x.$$

The Cartan relation is given by $\d G_n^1$;
according to lemma \ref{boundary}, \hfill\break
$\d G_{n+1}^1(x,x,y,y)=0$. Using relations (\ref{cellmn}), (\ref{defalpha}) and (\ref{umm}),
one has
\begin{align*}
0=&G^1_n\cdot (\Id+(2143)) (x,x,y,y)\\
+&\alpha_{n,n}(x,x,y,y)+\alpha_{n,0}\cdot (3214) (x,x,y,y)\\
=&2 G^1_n(x,x,y,y) + [u_n^1]_0 (x,x,y,y) \\
+&\sum_{s} [u_s^{s+1}]_{n-s} (y,x,x,y) +
[u_{s+1}^{s+1}]_{n-s-1} \cdot (3412) (y,x,x,y)\\
=&\sum_{l=0}^n e_0(e_{l},e_{n-l}\cdot (21)^l)(x,x,y,y)+
[u_0^1]_{n}(y,x,x,y) \\
+&\sum_{s=0}^{n-1} e_{n-s-1}\cdot(21)^{s}(x*y,y\cup_{s+1} x)
+e_{n-s-1}\cdot(21)^{s}(x\cup_{s+1} y,x*y) \\
\end{align*}
Using the commutativity of $*$ and $\cup_i$ for all $i$, one gets
\begin{align*}
\d G_{n+1}^1(x,x,y,y)=& \sum_{l=0}^n D_l(x)*D_{n-l}(y)+ D_n(x*y) \\
\end{align*}
which gives the Cartan relation (\ref{Cartanrel2}).

\medskip

The proof of Adem relation (\ref{Ademrel2}) relies on lemma \ref{lemmeuD}, and
on the relation $\d G^m_{n+1}(a,a,a,a)=0$. Combined with the relation
(\ref{cellmn}) one gets $\alpha_{n,p}(a,a,a,a)=\alpha_{n,n-p}(a,a,a,a)$ that is
\begin{multline*}
\sum_{s} [u^s_{p+s}(a,a,a,a)+u^{s+1}_{p+s}(a,a,a,a)]_{n-p-s}=\\
\sum_{t} [u^t_{n-p+t}(a,a,a,a)+u^{t+1}_{n-p+t}(a,a,a,a)]_{p-t}
\end{multline*}
\begin{multline*}
\Rightarrow \sum_s [(p-s,2s-1)+(p-s-2,2s+1)] D_{n-p-s} D_{\frac{p+s}{2}}(a)=\\
\sum_t [(n-p-t,2t-1)+(n-p-t-2,2t+1)] D_{p-t} D_{\frac{n-p+t}{2}}(a).
\end{multline*}
But $(x,y-2)+(x-2,y)=(x,y)$, hence
\begin{multline*}
\sum_s (p-s,2s+1) D_{n-p-s}D_{\frac{p+s}{2}}(a)=\\
\sum_t (n-p-t,2t+1) D_{p-t} D_{\frac{n-p+t}{2}}(a).
\end{multline*}
Since the first term is zero as soon as  $p-s$ is odd, we can
set $s=p-2l,$ and also $t=n-p-2k$. As a consequence,
\begin{multline*}
\sum_l (2l,2p-4l+1) D_{n-2p+2l}D_{p-l}(a)=\\
\sum_k (2k,2n-2p-4k+1) D_{2p-n+2k} D_{n-p-k}(a)
\end{multline*}
Using the 2-adic expansion, one gets
\begin{equation}\label{double}
(2l,2p-4l+1)=(l,p-2l);
\end{equation}
by setting $w:=p$ and $v:=n-p$, one obtains
\begin{align*}
\sum_l (l,w-2l) D_{v-w+2l}D_{w-l}(a)&=
\sum_k (k,v-2k) D_{w-v+2k} D_{v-k}(a)
\end{align*}
which is the Adem relation (\ref{Ademrel2}). ~~~\hfill{$\square$}

\subsubsection{Corollary}\label{structureB2}\it Let $\cal O$ be an
Adem-Cartan operad. The co\-ho\-mo\-lo\-gy of any $\cal O$-algebra is an
unstable level algebra over $\mathcal B_2$. \rm

\medskip

\noindent{\bf Proof.} By definition \ref{ADO}
 there is a morphism of operads $\LevAC\rightarrow \cal O$,
then it suffices to
prove the theorem for $\cal O=\LevAC$. Let $A$ be a
$\LevAC$-algebra, then $H^*(A)$ is a level algebra. In order to
prove the Adem-Cartan relations we compute the boundaries of
$G^1_n(a,a,b,b)$ and $G^m_n(a,a,a,a)$ for cocycles $a$ and $b$
that represents classes $[a]$ and $[b]$ in the cohomology. These
boundaries give the Adem-Cartan relations between $e_o(a,b)$,
$e_i(a,a)$, $e_j(b,b)$ which represent $[a]*[b]$, $D_i([a])$ and
$D_j([b])$ respectively.~~~~\hfill$\square$

\medskip

\noindent In particular for algebras over a $E_\infty$-operad, we recover the
results of May (\cite{May} and \cite{KM}).

\medskip

\subsubsection{Example} Given an $\mathcal E$-algebra, a
natural question is to know whether it is possible to extend this
structure into a structure of Adem-Cartan algebra. The following
example shows that it can not be done by imposing the triviality
of all $G^m_n$'s.
\\
Let us consider the torus $\mathbb T^2=S^1\times S^1$. Here is
the algebraic model of the normalized singular cochains
of $\mathbb T^2$ that we use:
$$A_{\mathbb T^2}=C^*(S^1)\otimes C^*(S^1).$$
The vector space $A_{\mathbb T_2}$ is generated by:
$1:=1\otimes 1$ in degree zero,
$\alpha:=a\otimes 1$ and $\beta:=1\otimes b$ in degree $1$
and $\alpha\beta=a\otimes b$ in degree 2.
The differential is trivial on $A_{\mathbb T^2}$.

\medskip

The $\E$-structure on $C^*(S^1)$ is given by $e_0(1,1)=1$, $e_1(a,a)=a$ where $a$ is the
generator in degree 1 and all the others are zero (see \cite{BF}). Hence the
 $\E$-structure on $A_{\mathbb T^2}$ is given by the coproduct $\psi$ on $\E$.
Since the differential is zero one has $\d (G_3^1(\alpha,\beta,\alpha,\beta))=0$.
But $\d G_3^1=G_2^1(\Id+(2143))+\alpha_{2,2}+\alpha_{2,0}\cdot (3214)$. Then, using the
definition of the $\alpha_{n,p}$'s and the commutativity of the $e_i$'s, one has
\begin{align*}
\d (G_3^1(\alpha,\beta,\alpha,\beta))=&G_2^1(\alpha,\beta,\alpha,\beta)+
G_2^1(\beta,\alpha,\beta,\alpha) \\
&+e_0(e_1,e_1\cdot(21))(\alpha,\beta,\alpha,\beta)+e_2(e_0,e_0)(\alpha,\beta,\alpha,\beta)\\
=&G_2^1(\alpha,\beta,\alpha,\beta)+
G_2^1(\beta,\alpha,\beta,\alpha) \\
&+e_1(\alpha,\beta)^2+e_2(\alpha\beta,\alpha\beta)
\end{align*}
\begin{align*}
\mbox{\rm But\ } e_1(\alpha,\beta)=&\psi(e_1)(a\otimes 1,1\otimes b)\\
=&(e_0\otimes e_1+e_1\otimes e_0\cdot(21)) ((a,1)\otimes(1,b))\\
=&e_0(a,1)\otimes e_1(1,b)+e_1(a,1)\otimes e_0(b,1)=0 \\
\mbox{\rm and\ } e_2(\alpha\beta,\alpha\beta)=&
(e_0\otimes e_2+e_1\otimes e_1\cdot (21)+e_2\otimes e_0)((a,a)\otimes(b,b)) \\
=&a\otimes b= \alpha\beta
\end{align*}
$$\mbox{\rm thus},\ \ \ \ G_2^1(\alpha,\beta,\alpha,\beta)+
G_2^1(\beta,\alpha,\beta,\alpha)=\alpha\beta$$
which proves that the action of $G_2^1$ is non zero.

\section{Operadic secondary cohomological operations.}

In this section we prove that there exist secondary
cohomological operations on the cohomology of an Adem-Cartan algebra $A$, and that these 
operations coincide with the Adams operations in case $A=C^*(X,\F2)$ is the singular cochain complex
of a space $X$.

\subsection{Secondary cohomological operations on Adem-Cartan algebras.} 
Let $\cal O$ be an Adem-Cartan operad and $A$ be an $\cal O$-algebra. Then $A$ is a 
$\LevAC$-algebra, hence $H^*(A)$ is endowed with a (non natural) structure of 
$\LevAC$-algebra (see corollary \ref{corinvariance}). Hence for any $(m,p)$ there
are morphisms 
$$\theta^{m,p}:H^n(A)\rightarrow H^{4n-p-1}(A)$$
 given by $\theta^{m,p}(x)=G_{p+1}^m(x,x,x,x)$.

\smallskip

Besides from corollary \ref{structureB2}, the cohomology $H^*(A)$ is an unstable level algebra over $\cal B_2$.
Let $x\in H^n(A)$ and
$$R_{Ad}(x)=\sum_i Sq^{m_i}Sq^{n_i}(x)=0$$
be an Adem relation with $x\in \underset{i}{\cap}\Ker(\Sq^{n_i})$.
Let $c\in A^n,\ \d c=0$ representing $x$. Then there exists
$(m,p)$ such that $(\d G_{p+1}^m)(c,c,c,c)=R_{Ad}(c)$:
$p=3n-m_i-n_i$ and $m$ depends on $n$ and $m_i+n_i$. Since
$\Sq^{n_i}(x)=0$, there exists $b_i\in A^{n+n_i-1}$ such that
$db_i=e_{n-n_i}(c,c)$. The element
$$b=\sum_{i}e_{n-m_i+n_i}(1,e_{n-n_i})(b_i,c,c)+e_{n-m_i+n_i-1}(b_i,b_i)$$
satisfies  $d(G^m_{p+1}(c,c,c,c)+b)=0$.

\subsubsection{Proposition} \label{indepbi} \it The class of $G^m_{p+1}(c,c,c,c)+b$
in \hfill\break $H^{m_i+n_i+n-1}(A)/\sum_i Im(\Sq^{m_i})$ does
not depend on the choices of the $b_i'$s and $c$. \rm

\medskip

\noindent {\bf Proof.} First it does not depend on the choices of
the $b_i$'s. Let $b'_i\in A^{n+n_i-1}$ such that $\d b'_i=\d
b_i=e_{n-n_i}(c,c)$, and
$$b'=\sum_{i}e_{n-m_i+n_i}(1,e_{n-n_i})(b'_i,c,c)+e_{n-m_i+n_i-1}(b'_i,b'_i)$$
Then $\d(b+b')=0$ and the following relation implies the result:
$$b+b'=\sum_i\Sq^{m_i}(b_i+b'_i)+ \d e_{n-m_i+n_i}(b_i+b'_i,b_i).$$

\smallskip

Secondly it does not depend on the choice of a representant $c$ of
$x$. Using the homotopy invariance principle as in
\ref{corinvariance}, there is a zig-zag of acyclic fibrations of
Adem-Cartan algebras
$$\xymatrix{
& \tilde A \ar[dl]_{f}\ar[dr]^{g}& \\
A   & & H^*(A)}$$ with $H^*(f) =H^*(g)$. Let $c$ and $c'$ be two
cocycles that represent $x$, and let $u,u'\in \tilde A$ such that
$\d u=0,\ f(u)=c,\ g(u)=x$ and the same for $u',c'$. Since
$\Sq^{n_i}(x)=0$, there exists $v_i,v'_i$ such that $\d
v_i=\Sq^{n_i}(u)$ and $\d v'_i=\Sq^{n_i}(u')$. Let $v,v'$ defined
as $b$. Let us prove that
$$[G^m_{p+1}(c,c,c,c)+f(v)]=[G^m_{p+1}(c',c',c',c')+f(v')]$$
in $H^{n+m_i+n_i-1}(A)/\underset{i}{\sum} Im(\Sq^{m_i}).$ This
is equivalent to prove that $g(v)+g(v')$ is in the sum of all
$Im(\Sq^{m_i})$. But
$$v=\sum_i e_{n-m_i+n_i}(v_i,\d v_i)+\Sq^{m_i}(v_i)$$
and $g(\d v_i)=0$ imply  that $g(v)=\sum_i \Sq^{m_i}(g(v_i)) \in
\sum Im(\Sq^{m_i})$, which proves the result.~~~~\hfill$\square$

As a consequence, we have defined a map
\begin{align*}
\psi^{m,p}:&\underset{i}{\cap}\Ker(\Sq^{n_i})\subset
H^{n}(A)&\longrightarrow&
{H^{n+m_i+n_i-1}(A)}/{\underset{i}{\sum} Im(\Sq^{m_i})}\\
&\quad\quad\quad\quad x&\mapsto&\quad\quad  [G^m_{p+1}(c,c,c,c)+b]
\end{align*}

\subsubsection{Proposition}\label{scogene}\it 
Let $\iota :\underset{i}{\cap}\Ker(\Sq^{n_i})\rightarrow
H^{n}(A)$ be the canonical inclusion and $\pi:H^{n+m_i+n_i-1}(A)\rightarrow
H^{n+m_i+n_i-1}(A)/\sum_i Im(\Sq^{m_i})$ be the canonical projection.
Then 
$$\pi\theta^{m,p}\iota=\psi^{m,p}.$$
\rm

\noindent{\bf Proof.} Using
the proof of the previous proposition, it only remains to show
that $[G_{p+1}^m(c,c,c,c)+b]=G_{p+1}^m(x,x,x,x)$ in
$H^{m_i+n_i+n-1}(A)/\sum_i Im(\Sq^{m_i})$ for $x\in
\underset{i}{\cap}\Ker(\Sq^{n_i})$. This is equivalent to prove
that $g(v)$ is in $\sum_i Im(\Sq^{m_i})$ for $v$ with $f(v)=b$,
which has been already proved.~~~~\hfill{$\square$}

\subsection{Adams operations.} Let $C^*(X;\mathbb F_2)$
denotes the singular cochains complex of a topological space $X$
and $H^*(X,\F2)$ its cohomology. We recall that $C^*(X;\mathbb
F_2)$ is an algebra over a $E_\infty$-operad \cite
{HS}, hence an Adem-Cartan algebra.

\medskip

Adams defined in an axiomatic way stable secondary cohomological
operations \cite {Ad}. His approach is topological, and uses the
theory of so-called ``universal examples''.
These operations correspond to Adem relations
$$R_{Ad}=\sum_{i}\Sq^{m_i}\Sq^{n_i},$$
and are denoted by $\Phi$. Let recall Adams axioms:
\\
{\bf Axiom 1.} For any $u\in H^n(X;\mathbb F_2)$, $\Phi(u)$ is
defined if and only if $\Sq^{n_i}(u)=0$ for all $n_i$.
\\
{\bf Axiom 2.} If $\Phi(u)$ is defined then
$$\Phi(u)\in H^{m_i+n_i+n-1}(X;\mathbb F_2)/\sum_i Im (\Sq^{m_i}).$$
\\
{\bf Axiom 3.} The operation $\Phi$ is natural.
\\
{\bf Axiom 4.} Let $(X,A)$ be a pair of topological spaces, we
have the long exact sequence
$$\ldots H^{n-1}(A;\mathbb F_2) \stackrel{\delta^*}{\rightarrow} H^n(X,A;\mathbb F_2)
\stackrel{j^*}{\rightarrow}H^n(X;\mathbb
F_2)\stackrel{i^*}{\rightarrow} H^n(A;\mathbb
F_2)\stackrel{\delta^*}{\rightarrow}\ldots$$ let $v\in
H^n(X,A;\mathbb F_2)$ be a class such that $\phi$ is defined on
$j^*(v)\in H^n(X;\mathbb F_2)$. Let $w_i\in H^*(A;\mathbb F_2)$
such that $\delta^*(w_i)=\Sq^{n_i}(v)$. Then, we have
$$i^*\Phi(j^*(v))=\sum_i\Sq^{m_i}(w_i)\in H^*(A;\mathbb F_2)/i^*(\sum_i Im (\Sq^{m_i})).$$
\\
{\bf Axiom 5.} The operation $\Phi$ commutes with suspension.

\medskip

Later on, Kristensen proved that these operations
can be defined at the cochain level, using the existence of a
coboundary which creates the stable secondary cohomological
operation defined by Adams (\cite {Kr}, chapter $6$). More
precisely, for an Adem relation $R_{Ad}$ and a class $x\in
\cap_i\Sq^{n_i}$, Kristensen defines cochain operations $\theta$
such that the differential of $\theta(c)$ (c is a representant of
$x$) gives a cocycle representing an Adem relation  $R_{Ad}(x)$.
If one chooses $b_i$ such that $db_i=\Sq^{n_i}(c)$, then one gets
a cocycle, and a cohomology class
$$Qu^r(c)=[\theta(c)+\sum_i(e_{n-m_i+n_i}(1,e_{n-n_i})(b_i,c,c)+e_{n-m_i+n_i-1}(b_i,b_i))].$$
Then,

\subsubsection{\bf Theorem}(Kristensen, theorem 6.1 of \cite{Kr})\label{thKr}
{\it Any operation $x\mapsto Qu^r(c)$ satisfies axiom 1-5 of Adams.}

\medskip

\subsubsection{Corollary}\label{AdamsId}\it 
The maps $\psi^{m,p}$ coincide with the stable secondary
cohomological operations of Adams. \rm

\medskip

\noindent {\bf Proof.} The proof relies on theorem \ref{thKr} with
$\theta(c)=G^m_{p+1}(c,c,c,c)$.~~~~\hfill{$\square$}

\subsubsection{Theorem}\label{cohop2} \it The stable secondary cohomological
operations $\phi$ of Adams extend to maps $\theta^{m,p}:
H^n(X,\mathbb F_2)\rightarrow H^{n+m_i+n_i-1}(X,\F2)$. More
precisely, if we denote by
$\iota :\underset{i}{\cap}\Ker(\Sq^{n_i})\rightarrow H^n(X,\F2)$ and
$\pi:H^{n+m_i+n_i-1}(X,\F2)\rightarrow
H^{m_i+n_i+n-1}(X,\F2)/\sum_i Im(\Sq^{m_i})$ then
$$\phi=\pi\theta^{m,p}\ \iota.$$ \rm

\noindent{\bf Proof.} It is the translation of theorem \ref{scogene} 
for $A=C^*(X,\F2)$.~~~~\hfill{$\square$}

\section{Proof of technical lemmas} \label{technique}

In this section proposition
\ref{suite} and lemma \ref{lemmeuD} are proved.

\subsection{Lemma} \label{lemmeij} \it For any $i,j \leq 2^p-1$ one has
\begin{align*}
(i,j)=& 0,\  ~{\rm if}~ i+j \geq 2^p,\ ~{\rm and}~ \\
(i,j)=& (2^p-i-j-1,j).
\end{align*}
\rm

\noindent{\bf Proof.} Let $\sum_{l=0}^{p-1} a_l2^l$ and $\sum_{l=0}^{p-1}
b_l2^l$  be the 2-adic expansion of $i$ and $j$ respectively.
Recall that $(i,j)=1$ if and only if the 2-adic expansion of $i+j$ is
$\sum(a_l+b_l)2^l$. If $i+j\geq 2^p$ this is not the case, thus $(i,j)=0$.
\\
If $(i,j)=1$ then the 2-adic expansion of $2^p-1-i-j$ is
$\sum_{l=0}^{p-1}(1-a_l-b_l)2^l$,
thus the 2-adic expansion
of $(2^p-1-i-j)+j$ has for coefficients $(1-a_l-b_l) +(b_l)$.
Consequently $(2^p-1-i-j,j)=1$. The converse is true
by symmetry.
\\
Note that the first assertion is a consequence of the second one,
because if $(i+j)\geq 2^p$ then $(2^p-1-i-j)<0$, and
$(\alpha,\beta)=0$ if $\alpha<0$ or $\beta<0$.
~~~~~~\hfill{$\square$}

\subsection{Proof of lemma \ref{lemmeuD}.}
Using the commutativity of $\cup_x$ one has
\begin{multline*}
[u^m_n]_x(a,a,a,a)= \sum_{i=0}^{2^{k+1}-1} \sum_{ 0\leq 2^{k+1}\delta-i\leq n}\\
(n-m+i,m-1)(i,m)
D_{2^{k+1}\delta-i}(a)\cup_x D_{n+i-2^{k+1}\delta}(a)\\
+(n-m+i,m-1)(i-1,m) D_{2^{k+1}\delta-i}(a)\cup_x D_{n+i-2^{k+1}\delta}(a)\\
= \sum_{i,\delta}
(n-m+i,m-1)(i,m-1) D_{2^{k+1}\delta-i}(a)\cup_x D_{n+i-2^{k+1}\delta}(a).\\
\end{multline*}
Let $0\leq j\leq 2^{k+1}-1$ such that
$n+i\equiv -j \ [2^{k+1}]$, then there exists $\delta'$ such that $n+i-2^{k+1}\delta=2^{k+1}\delta'-j$,
and lemma \ref{lemmeij} implies
\begin{equation}
\begin{split}\label{gauchedroite}
(n-m+i,m-1)&=(2^{k+1}(\delta+\delta')-m-j,m-1)\\
&=(j,m-1)\\
(n-m+j,m-1)&=(i,m-1)
\end{split}
\end{equation}
As a consequence, if $i\not=j$ or $i=j$ and $\delta\not=\delta'$,
the 2 following terms in $[u_n^m]_x(a,a,a,a)$
\begin{multline*}
(n-m+i,m-1)(i,m-1) D_{2^{k+1}\delta-i}(a)\cup_x D_{2^{k+1}\delta'-j}(a)+\\
(n-m+j,m-1)(j,m-1) D_{2^{k+1}\delta'-j}(a)\cup_ x D_{2^{k+1}\delta-i}(a)
\end{multline*}
vanish. Hence if there exists $(i,\delta)$ such that
$2^{k+1}\delta-i=n-(2^{k+1}\delta-i)$, then
$$[u^m_n]_x(a,a,a,a)=(n-m+i,m-1)^2 D_xD_{\frac{n}{2}}(a).$$
Relation (\ref{double}) implies
$(n-m+i,m-1)=(2n-2m+2i,2m-1)=(n-2m,2m-1).$ Furthermore, if
$(n-2m,2m-1)=1$, then $n$ is even and we can pick $0\leq r\leq 2^{k+1}-1$ such that
$r\equiv -\frac{n}{2} \ [2^{k+1}]$.~~~~~~\hfill
{$\square$}

\subsection{Proof of proposition \ref{suite}} We have to prove the relations
(\ref{unmzero}), (\ref{umm}), (\ref{du}) and (\ref{dalpha}). The relation (\ref{dalpha}) is
straightforward using relation (\ref{du}) and definition (\ref{defalpha}).

\medskip

{\sl Proof of relations (\ref{unmzero}) and (\ref{umm})--} Assume $n\leq m$.
The condition $0\leq 2^{k+1}\delta -i \leq n\leq m \leq 2^{k+1}-1$
implies $\delta$ equals $0$ or $1$. If $\delta=0$, then $i=0$.
If $\delta=1$, then $i+m\geq 2^{k+1}$ and $(i,m)=0$ by lemma \ref{lemmeij}. So $[u^m_n]_x$ writes
\begin{multline*}
[u_n^m]_x=(n-m,m-1)e_x\cdot(21)^{m-1}(e_0,e_n) \\
+\sum_{i=0}^{2^{k+1}-1}(n-m+i,m-1)(i-1,m)e_x\cdot(21)^{m-1}(e_{2^{k+1}-i},e_{n+i-2^{k+1}}\cdot(21))
\end{multline*}
If $n<m$, then $(n-m,m-1)=0$ and $i-1+m\geq 2^{k+1}$ implies $(i-1,m)=0$,which proves relation (\ref{unmzero}).

If $n=m$, $(i-1,m)=0$ for all $i\not=2^{k+1}-m$ and
$(2^{k+1}-m,m-1)(2^{k+1}-1-m,m)=(0,m-1)(0,m)$ by virtue of lemma \ref{lemmeij}.
This proves relation (\ref{umm}).

\medskip

{\sl Proof of relation (\ref{du})--} For the convenience of the reader, let
$$B_{x,[m],n}^{\delta,i}=e_x\cdot(21)^{m-1}
(e_{2^{k+1}\delta-i},e_{n+i-2^{k+1}\delta})$$
where $[m]$ means $m$ mod 2, then
$$[u_n^m]_x=\sum_{i,\delta}(n-m+i,m-1)[ (i,m)B_{x,[m],n}^{\delta,i} + (i-1,m)B_{x,[m],n}^{\delta,i}\cdot(1243)].$$

{\sl Remarks:}
\begin{itemize}
\item[a)] Let $P=\{\Id,(2134),(2143),(1243)\}\subset \Sigma_4$. Then the
set $\mathcal F=\{B_{x,[m],n}^{\delta,i}\cdot \sigma,
B_{x',[m+1],n'}^{\delta',j}\cdot \tau,\forall
x,n,x',n',i,i',\delta,\delta',
\forall \sigma,\tau \in P\}$
form a free system in $\E(4)$.
\item[b)] For any $\sigma\in P$,
$B_{x,[m],n}^{\delta,i}\cdot\sigma F=B_{x,[m],n}^{\delta,i}\cdot
F\sigma=B_{x,[m+1],n}^{\delta,i}\cdot\sigma$.
\end{itemize}

There are two cases to consider: if
$2^{k}\leq m \leq 2^{k+1}-2$ (then $m+1\leq 2^{k+1}-1$) or if  $m=2^{k+1}-1$. Since computation
are long but not difficult,
we'll present only the first case.

\begin{align*}
\d [u_{n+1}^{m+1}]_{x+1}=&\sum_{i,\delta}
\underbrace{(n-m+i,m)(i,m+1)}_{a_{i}}[B_{x+1,[m+1],n}^{\delta,i+1}\cdot(\Id+(2134)) \\
+&B_{x+1,[m+1],n}^{\delta,i}\cdot(\Id+(1243))+B_{x,[m+1],n+1}^{\delta,i}\cdot(\Id+F)]\\
+&\underbrace{(n-m+i,m)(i-1,m+1)}_{b_{i}}[B_{x+1,[m+1],n}^{\delta,i+1}\cdot((1243)+(2143)) \\
+&B_{x+1,[m+1],n}^{\delta,i}\cdot(\Id+(1243))+B_{x,[m+1],n+1}^{\delta,i}\cdot(\Id+F)(1243)]\\
\end{align*}
\begin{multline*}
[u_n^{m+1}]_{x+1}(\Id+((2143))=\\
\sum_{i,\delta}\underbrace{(n-m+i-1,m)(i,m+1)}_{c_i} B_{x+1,[m+1],n}^{\delta,i}(\Id+(2143))\\
+\underbrace{(n-m+i-1,m)(i-1,m+1)}_{d_i} B_{x+1,[m+1],n}^{\delta,i}((1243)+(2134))\\
\end{multline*}
\begin{multline*}
[u_n^{m}]_{x+1}(F+(4321))=\\
\sum_{i,\delta}
\underbrace{(n-m+i,m-1)(i,m)}_{e_i} B_{x+1,[m],n}^{\delta,i}(F+(4321))\\
+\underbrace{(n-m+i,m-1)(i-1,m)}_{f_i} B_{x+1,[m],n}^{\delta,i}(1243)(F+(4321))\\
\end{multline*}
Note that
$$B_{x+1,[m],n+1}^{\delta,i}\cdot(4321)=
e_{x+1}\cdot (21)^{m}(e_{n+i-2^{k+1}\delta},e_{2^{k+1}\delta-i})
\cdot(3412)(4321)$$
Hence by using relation (\ref{gauchedroite}) we get
\begin{multline*}
\sum_{i,\delta}(n-m+i,m-1)(i,m) B_{x+1,[m],n}^{\delta,i}\cdot(4321)=\\
\sum_{j,\delta'}\underbrace{(n-m+j-1,m)(j,m-1)}_{m_j}B_{x+1,[m+1],n}^{\delta',j}(2143)
\end{multline*}
and
\begin{multline*}
\sum_{i,\delta}(n-m+i,m-1)(i-1,m) B_{x+1,[m],n}^{\delta,i}(1243)(4321)=\\
\sum_{j,\delta'}\underbrace{(n-m+j,m)(j,m-1)}_{l_j}B_{x+1,[m+1],n}^{\delta',j}(1243)
\end{multline*}
\begin{multline*}
[u_{n+1}^{m+1}]_x(1+F)=\sum_{i,\delta}
\underbrace{(n-m+i,m)(i,m+1)}_{g_i} B_{x,[m+1],n+1}^{\delta,i}(\Id+F)\\
+\underbrace{(n-m+i,m)(i-1,m+1)}_{h_i}
B_{x,[m+1],n+1}^{\delta,i}(1243)(\Id+F)\\
\end{multline*}
Thus to prove relation (\ref{du}), it suffices to prove that the sum
of all the coefficients of elements of ${\mathcal F}$ vanishes.
For instance, the coefficient of $B_{x+1,[m+1],n}^{\delta,i}\cdot\Id$
is $a_{i-1}+a_i+b_i+c_i+e_i$ that is
\begin{align*}
&(n-m+i-1,m)(i-1,m+1)+(n-m+i,m)(i,m+1)+\\
&(n-m+i,m)(i-1,m+1)+(n-m+i-1,m)(i,m+1)+\\
&(n-m+i,m-1)(i,m)\\
=&(n-m+i,m-1)(i-1,m+1)+(n-m+i,m-1)(i,m+1)+\\
&(n-m+i,m-1)(i,m)=0.\\
\end{align*}
The coefficient of
$B_{x+1,[m+1],n}^{\delta,i}\cdot (2143)$
is  $b_{i-1}+c_i+m_i$ that is
\begin{align*}
&(n-m+i-1,m)(i-2,m+1)+(n-m+i-1,m)(i,m+1)\\
&+(n-m+i-1,m)(i,m-1)=0.
\end{align*}
All the others computation go the same. ~~~~\hfill{$\square$}

\medskip
{\bf Acknowledgements.} The authors would like to thank Benoit
Fres\-se for fruitful discussions. They are especially indebted
to Lionel Schwartz for pointing them the importance of level
algebras in algebraic topology. The first author would like to
thank Marc Aubry, Yves Felix and Jean-Claude Thomas for their
support, the paper was written during his stay at the University
of Louvain-La-Neuve and at the CRM (Barcelona). The second author
would like to thank the University of Louvain-La-Neuve and the CRM
for their hosting.

\bigskip

\end{document}